\magnification=1200

\hsize 17truecm
\vsize 23truecm

\font\twelvec=msbm10 at 10pt
\font\sevenc=msbm10 at 7pt
\font\fivec=msbm10 at 5pt

\newfam\co
\textfont\co=\twelvec
\scriptfont\co=\sevenc
\scriptscriptfont\co=\fivec

\def\const{\mathop{\rm const.}\nolimits}

\def\const{\mathop{\rm const.}\nolimits}
\def\det{\mathop{\rm det}\nolimits}

\def\ess{\mathop{\rm ess}\nolimits}
\def\exp{\mathop{\rm exp}\nolimits}

\def\Id{\mathop{\rm Id}\nolimits}

\def\id{\mathop{\rm Id}\nolimits}
\def\im{\mathop{\rm Im}\nolimits}

\def\e{\mathop{\rm e}\nolimits}

\def\ker{\mathop{\rm Ker}\nolimits}
\def\lim{\mathop{\rm lim}\nolimits}

\def\mod{\mathop{\rm mod}\nolimits}

\def\re{\mathop{\rm Re}\nolimits}

\def\supp{\mathop{\rm supp}\nolimits}
\def\Tr{\mathop{\rm Tr}\nolimits}

\def\inf{\mathop{\rm inf}\nolimits}

\def\Op{\mathop{\rm Op}\nolimits}

\def\WF{\mathop{\rm WF}\nolimits}

\def\Sum{\displaystyle\sum}
\def\e{\mathop{\rm \varepsilon}\nolimits}

\baselineskip 15pt

\centerline{\bf SEMI-CLASSICAL QUANTUM MAPS OF SEMI-HYPERBOLIC TYPE}
\bigskip
\centerline{Hanen LOUATI ${}^{1}$, Michel ROULEUX ${}^{2}$}
\medskip
\centerline {${}^{1}$ Universit\'e de Tunis El-Manar, D\'epartement de
Math\'ematiques, 1091 Tunis, Tunisia}

\centerline {e-mail: louatihanen42@yahoo.fr}

\centerline {${}^{2}$ Aix Marseille Universit\'e, Universit\'e de Toulon, CNRS, CPT,
Marseille, France}

\centerline {e-mail: rouleux@univ-tln.fr}
\bigskip
\noindent {\it Abstract}: 
Let $M={\bf R}^n$ or possibly a Riemannian, non compact manifold. 
We consider semi-excited resonances for a $h$-differential operator 
$H(x,hD_x;h)$ on $L^2(M)$ 
induced by a non-degenerate periodic orbit $\gamma_0$ of semi-hyperbolic type,
which is contained in the non critical energy surface $\{H_0=0\}$.
By semi-hyperbolic, we mean that the linearized Poincar\'e map $d{\cal P}_0$
associated with $\gamma_0$ has at least one eigenvalue of modulus greater (or less) than 1,
and one eigenvalue of modulus equal to 1,
and by non-degenerate that 1 is not an eigenvalue, which implies a family $\gamma(E)$ with the same properties.
It is known that an infinite number of periodic orbits
generally cluster near $\gamma_0$, with periods approximately multiples of its primitive period. 
We construct the monodromy and Grushin operator, adapting some arguments by [NoSjZw], [SjZw], and compare 
with those obtained in [LouRo], which ignore the additional orbits near $\gamma_0$, but still give the right 
quantization rule for the family $\gamma(E)$. 

\medskip
\noindent {\bf 1. Introduction}
\smallskip
Let $M$ be a smooth manifold (for simplicity here $M={\bf R}^n$, but our results hold in more general cases, see Examples 1
and 2 below), and $H(y,hD_y;h)$ be a semi-classical Partial Differential Operator of second order, we assume to be self-adjoint  on $L^2(M)$, and 
satisfy usual hypotheses required in the framework of resonances. In particular, 
its Weyl symbol $H(y,\eta;h)$, in the sense of $h$-$\Psi$DO, belongs to the class
$$S^N(\langle\eta\rangle^2)=\{H\in C^\infty(T^*{\bf R}^n): \forall\alpha\in{\bf N}^{2n}, \exists C_\alpha>0,
|\partial^\alpha_{(y,\eta)}H(y,\eta;h)|\leq C_\alpha h^N \langle\eta\rangle^2 \}$$
i.e. is of growth at most quadratic in momentum at infinity (here $\langle\eta\rangle^2=1+|\eta|^2$). 
\smallskip
\noindent {\it a) Main hypotheses}
\smallskip
\noindent $\bullet$ {\bf Hypothesis 1} (Ellipticity, regularity of coefficients and behavior at infinity).

$H$ is elliptic (i.e. $|H(y,\eta;h)+i|\geq \const \langle\eta\rangle^2$)
and extends analytically in a ``conic'' neighborhood of the real domain
$$\Gamma_0=\{(y,\eta) \in T^*{\bf C}^n : |\im (y,\eta)|\leq\const \langle\re (y,\eta)\rangle\}\leqno(1.1)$$ 
where it has the semi-classical expansion
$$H(y,\eta; h) \sim H_0(y,\eta)+hH_1(y,\eta)+\cdots, h\rightarrow 0$$
To fix the ideas, we assume $H(y,\eta;h)-\eta^2\to 0$ in $\Gamma_0$ when $|\re y|\to\infty, |\im y|\leq\const \langle\re y\rangle$.
This assumption can be relaxed, see [HeSj].

Then $\inf\sigma_{\ess }(H(y,hD_h;h))=0$ (actually $H$ has only continuous spectrum above 0)
and we define the {\it resonances} of $H$ near $E_0>0$ by the method of analytic distorsions, as the discrete spectrum
of some non self-adjoint extension of $H$.  

Namely, let $\Gamma\subset{\bf C}^n$ be a real, totally real
submanifod of dimension $n$, and  $H(y,hD_y;h)=\Sum_{|\alpha|\leq2}a_\alpha(z;h)(hD_z)^\alpha$ a differential operator 
with $C^\infty$ coefficients in some suitable complex neighborhood $\Gamma^{\bf C}$ of $\Gamma$ (here $D_z$ denote the holomorphic
derivative with respect to coordinates in $\Gamma^{\bf C}$). Then we can define a differential operator $H_\Gamma:C^\infty(\Gamma)
\to C^\infty(\Gamma)$, such that, if $u$ is holomorphic, then $(Hu)|_{\Gamma}=H_\Gamma(u|_\Gamma)$. Now assume that $H(y,\eta;h)$ 
is defined in $\Gamma_0$ as in (1.1). For $0\leq\theta\leq\theta_0$, we let $\Gamma=\Gamma_\theta$ be parametrized by $f_\theta\in
C^\infty({\bf R}^n;{\bf C}^n)$ such that $f_\theta(y)=y$ for $y$ in a compact set and $f_\theta(y)=e^{i\theta}y$ for large $y$.
The corresponding family of operators $H_\theta=H_{\Gamma_\theta}$ on $L^2$ is known to be an analytic family of type (A) and
$\sigma_{\ess }(H_\theta)=e^{-2i\theta}{\bf R}^+$. Moreover, when $\theta>0$, $H_\theta$ is Fredholm and may also have discrete
eigenvalues in the lower-half plane near $E_0$, called (outgoing) resonances. The resonant (or extended) states are the
associated eigenfunctions. See [Co], [Va], [ReSi], [BrCoDu], [HeSj], [G\'eSi] for related approaches, which turn out
to be essentially equivalent ([HeMa]). We follow here mainly [NoSjZw].

Since we shall mostly consider $H(y,hD_y;h)$ as a $h$-$\Psi$DO, we shall rather denote it by $H^w(y,hD_y;h)$. 

Locating precisely resonances near $E_0$ (like Bohr-Sommerfeld quantization conditions)
hinges on properties of the Hamiltonian flow on the energy surfaces nearby $E_0$. As recalled briefly in Appendix, 
we need to choose 
distorsion $f_\theta$ accurately, as well as other phase-space distorsions, or Lagrangian deformations.

\smallskip
\noindent $\bullet$ {\bf Hypothesis 2} (Regularity of energy surface)

To save notations we change $H_0$ to $H$ when considering
classical quantities. 
We fix a regular energy surface $\{ H(y,\eta)=0\}$, and assume there is an energy interval $I$ around $E_0$, so that the Hamilton vector field 
$X_{H}$ has no fixed point on $\{ H(y,\eta)=E\}$, for $E\in I$.  

Let 
$\Phi^t =\exp(tX_H) : T^* {\bf R}^n \rightarrow T^*{\bf R}^n$ be the Hamiltonian flow and 
$${\cal K}(E) = \{\rho \in T^\ast {\bf R}^n, H_0(\rho)=E, \Phi^{t}(\rho) \
\hbox{doesn't grow to infinity} \ \hbox{as} \ |t|\to\infty\}\leqno(1.2)$$
the trapped set at energy $E$. Simplest situation holds when ${\cal K}(E)$ is a fixed point [BrCoDu]. 

\smallskip
\noindent $\bullet$ {\bf Hypothesis 3}  (Trapped set at energy 0) 

We assume here that ${\cal K}_0={\cal K}(0)$ contains a periodic orbit of primitive period $T_0$. 
The differential of Poincar\'e map (or first return map) is a
symplectic automorphism ${\cal P}_0$ of the normal space $\Sigma_0=\Sigma(0)$ of $\gamma_0$ in $H^{-1}(0)$,
is a manifold of dimension $2d=2(n-1)$, which is called Poincar\'e section. 

Let $\lambda\in{\bf C}$ be the eigenvalues of $A=d{\cal P}_0$ (or Floquet multipliers). 
The periodic orbit $\gamma_0$ is said 
non degenerate if 1 is not a Floquet multiplier. 
By Poincar\'e Continuation Theorem,
there is a one parameter family of periodic orbits $\gamma(E)\subset H_0^{-1}(E)$ containing $\gamma_0$, and $\gamma(E)$ is non-degenerate
for $E$ small enough. By abuse of notations, we shall still call Poincar\'e section the smooth 
foliation $\overline{\Sigma}=\bigcup_{E\in I} \Sigma(E)$ transverse to
$\overline{\gamma}=\bigcup_{E\in I}\gamma(E)$. 

An eigenvalue $\lambda\in {\bf C}$ of $A$ is called elliptic (ee) if
$|\lambda|=1$ ($\lambda \neq \pm 1$) and hyperbolic (he) if $|\lambda| \neq 1$. The corresponding eigenspace will be denoted by 
$F_\lambda$. 
\smallskip
\noindent $\bullet$ {\bf Hypothesis 4} (genericity properties of linearized Poincar\'e map)
$$F_{\pm1} = \{0\}, \quad F_{\lambda} = \{0\}, \ \forall \lambda \leq 0$$
 
In particular, 
we can define $B = \log A$. Eigenvalues 
$\mu=\mu(\lambda)$ of $B$ (Floquet exponents) verify $\mu(\overline{\lambda})=\overline{\mu(\lambda)}$, 
$\mu=\log\lambda$. Exponent $\mu$ is said ee if $\re\mu=0$, 
real-hyperbolic (hr) if $\im\mu=0$, loxodromic or complex-hyperbolic (hc) if $\re\mu\neq0, \im\mu\neq0$.

Eigenvalues of $B$ have the form $\mu_j, -\mu_j,\overline{\mu_j},-\overline{\mu_j}\neq 0$, $\re \mu_j\geq0$,
with same multiplicity. 
For simplicity, assume eigenvalues $\mu_j$ are distinct.
\smallskip
\noindent $\bullet$ {\bf Hypothesis 5} (Hyperbolicity)

We are interested in the case where $\gamma_0$ is unstable:
$A$ is hyperbolic, i.e. has at last one eigenvalue $|\lambda|\neq1$.
We say we have pure (or complete) hyperbolicity iff $\re\mu_j > 0$ for all $j$. 
In case of complete hyperbolicity, $\gamma_0$ is isolated, and we will assume (excluding e.g. symmetries, which would involve tunneling,
as in Example 2 below)
that the trapped set reduces to $K_0=\gamma_0$. 
We say we have partial, or semi-hyperbolicity, iff there exists both $j$ with $\re\mu_j >0$
and $k$ with $\re\mu_k=0$. This is generically the case for hyperbolic systems [Ar]. 

Recall the Center/Stable/Unstable manifolds Theorem (see e.g. [G\'eSj2] for a review):
Let $\gamma_0$ be a non-degenerate periodic orbit, and $\gamma(E)$
as above.  Then there exists a closed symplectic submanifold $C$ ({\it center manifold}) containing $\gamma(E)$ 
and such that $X_{H_0}(\rho)$ is tangent to $C$ at every point $\rho\in C$ 
(i.e. $C$ invariant under the
Hamiltonian flow). There exist also two vector bundles 
$N^\pm$, such that for all $\rho\in\Sigma$, $T_\rho C^\sigma=N^+_\rho\oplus N^-_\rho$
(here $T_\rho C^\sigma$ is the orthogonal symplectic of $T_\rho C$). Moreover, $N^\pm|_{\gamma(E)}$ are invariant under the
Hamiltonian flow, which is contracting on $N^-|_{\gamma_E}$ and expansive on $N^+|_{\gamma_E}$.
Note that hr and hc components which belong to outgoing/incoming manifolds,
differ only by technical aspects, while ee components which belong to the center manifold, play a distinct role.
 
Here we say that $\gamma_0$ (and hence the family $\gamma(E)$) is {\it unstable} (e.g. in Lyapunov sense) if $N^\pm\neq0$, i.e.
when $A$ is hyperbolic. 

\smallskip
\noindent $\bullet$ {\bf Hypothesis 6} (Strong non-resonance condition and twist condition) 

$$\forall k\in {\bf Z}^d: \ \sum_{j=1}^{d} k_j\mu_j\in 2i\pi{\bf Z}  \Longrightarrow  k_j = 0, \forall j$$

Let $0<\ell<d$ be the number of elliptic elements, i.e. $\mu_j\in i{\bf R}, \im\mu_j>0$. 
Assume moreover that ${\cal P}_0|_C$ ($C$ the center manifold) is of twist type, i.e.  
the non linear Birkhoff invariants, are non degenerate. In particular $\gamma(E)$ is $N$-fold non- degenerate for all $N$. 
By Lewis-Birkhoff Fixed Point Theorem, see [Kl,Thm.3.3.3], in every neighborhood of $(0,0)\in{\bf R}^{2\ell}$, there exists infinitely many 
periodic points (i.e. belonging to periodic orbits). The number of orbits of bounded period is finite. 

Applying this Theorem to the normally hyperbolic symplectic invariant manifold $C$ for Poincar\'e map,
we find a sequence of 
periodic orbits $\gamma_k$ with (primitive) periods $T_k\to\infty$ clustering on $\gamma_0$, with $T_k\approx kT_0$
in the limit $k\to\infty$, and $\supp(\gamma_k)\to\supp(\gamma_0)$
(``infra-red limit'') So we assume the trapped set is of the form (excluding, as we already pointed out, 
other components of ${\cal K}_0$ in $\{H=0\}$)
$${\cal K}_0=\overline{\bigcup_{k\in{\bf N}}\gamma_k}\leqno(1.5)$$
In particular ${\cal K}_0$ is topologically 1-D, and only one Poincar\'e section (or 2 equivalent Poincar\'e sections)
is needed to describe the dynamics near ${\cal K}_0$, which simplifies the situation presented in [NoSjZw].
Moreover, there is structural stability of ${\cal K}(E)$ [KaHa,Thm.18.2.3]: namely the flows $\Phi^t|_{{\cal K}(E)}$ and 
$\Phi^t|_{{\cal K}_0}$ are conjugated, up to time reparametrization, by a homeomorphism close to identity. For instance it could happen that
${\cal K}(E)=\overline{\bigcup_{k\in{\bf N}}\gamma_k(E)}$, but this is not actually needed, for orbits with large period are unstable.  
It follows that we can choose Poincar\'e section $\overline\Sigma$ transverse to
$\overline{\cal K}=\bigcup_{E\in I}{\cal K}(E)$. We shall assume, as in [NoSjZw], that $\partial\overline\Sigma$ does not intersect 
$\overline{\cal K}$.
Our situation is very similar to [NoSjZw], and the more simple structure of the flow allows for some simplifications of the proof.
\medskip
\noindent{\it b) Examples} 
\smallskip
1) Poincar\'e example of a pure hyperbolic orbit: $H(y,hD_y)$ is the geodesic flow on one-sheeted hyperboloid in ${\bf R}^3$ (``diabolo''): 
the throat circle $\gamma_0$ is an unstable hr periodic orbit (geodesic).
\smallskip
2) The geodesic flow $H(y,hD_y)$ on a surface of revolution $M$ embedded in ${\bf R}^4$ with axis $Oz$, projecting on $x,z$-variables 
as the ``double diabolo'', a surface homeomorphic to the one-sheeted hyperboloid in ${\bf R}^3$, but with two throat circles, separated by 
a crest circle.  The effective Hamiltonian has principal symbol $H_0=\eta^2+\zeta^2+\bigl((2z^4-z^2+1)\cosh y\bigr)^{-2}-1$,
with $(x,\xi)$ as cyclic variables. For some energy $E_1$ there are two periodic geodesics of hyperbolic type 
(the throat circles) situated symmetrically on the hyperplanes $z=\pm{1\over4}$ (with two hr pairs); 
our constructions apply modulo tunneling corrections. 
For some energy $E_2<E_1$ there is  
one periodic geodesic of semi-hyperbolic type (the crest circle), with one ee pair and one hr pair. 
This is the generic situation, unlike the purely hyperbolic case as in previous Example. See [Chr2,App.C]. 
\smallskip
3) $H(y,hD_y)=-h^2\Delta_y+|y|^{-1}+ay_1$ on $L^2({\bf R}^n)$ (repulsive Coulomb potential perturbed by Stark effect) near an energy level
$E>2/\sqrt a$, or more generally, Schr\"odinger operators with potentials with two or more bumps. 
Their periodic orbits are generally hr (also called {\it librations}). See [G\'eSj], [Sj3].
\smallskip
4) Non-autonomous case [Tip]: Atom in a periodically polarized electric field $H_1(t)=-h^2\Delta+V(|x|)+E(t)\cdot x$ on $L^2({\bf R}^3)$,
$E(t)=\cos\omega t\widehat x_1+\sin \omega t\widehat x_2$. 
After some transformation, Floquet operator takes the form
$$U(s+T,s)=e^{-i\omega (s+T) L_{x_3}/h}e^{iTP(x,hD_x)/h}e^{i\omega (s+T) L_{x_3}/h}$$
where $T={2\pi\over\omega}$, 
$p(x,\xi)=(\xi-a)^2+V(|x|)-\omega(x_1\xi_2-x_2\xi_1)$, and 
$a=(1/\omega,0,0)$. The operator $e^{iTP(x,hD_x)/h}$ is now independent of time, and plays the role of the monodromy operator constructed below.
\medskip
\noindent {\it c)  Main result on resonances in the semi-hyperbolic case}
\smallskip
Our main result for which we sketch a proof in Sect.2,  is a straighforward generalization of
[NoSjZw], when allowing for elliptic Floquet exponents, and of [G\'eSj1] in the hyperbolic case. We summarize it as follows.
\medskip
\noindent {\bf Theorem 1.1}:  Under the Hypotheses 1-6 above, consider the spectral window
${\cal W}_h=[E_0-\e _0,E_0+\e _0]-i]0,Ch\log(1/h)]$. Then if $\e _0,C>0$ are small enough, there is $h_0>0$ small enough
and a family of matrices $N(z,h)$, such that  
the zeroes
of $\det(\id-N(z,h))$ give all resonances
of $H^w(x,hD_x;h)$ in
${\cal W}_h$ with correct multiplicities. 
The matrices $N(z;h)$ of order $k_h\sim h^{-n+1}$ are of the form $N(z;h)=\Pi_h{\cal M}(z;h)\Pi_h+{\cal O}(h^N)$
where $\Pi_h:L^2\to L^2$ (the weighted Hilbert space) are projectors of rank $k_h$ and ${\cal M}(z;h)$ is the monodromy 
operator quantizing Poincar\'e map and computed in Sect.2 below. 

\bigskip
\noindent{\it d) Bohr-Sommerfeld (BS) quantization rules}
\smallskip
BS for an hyperbolic orbit $\gamma(E)$ are known for a long time, see [G\'eSj1], [Vo]; in [LouRo1,2] we use the method
presented in Sect.3 below, ignoring the orbits accumulating on $\gamma(E)$. Our proof holds {\it stricto sensu} only
in the complete hyperbolic case, but the result turns out to be correct otherwise,
provided we consider only resonances associated with the family $\gamma(E)$. A peculiarity of BS rules for resonant spectrum
is that they cannot be simply derived from the construction of quasi-modes as in the self-adjoint case (see e.g. [BLaz]). We have:
\smallskip
\noindent{\bf Theorem 1.2} [LouRo]: 
Under the hypotheses above, let us define the {\rm semi-classical} action along $\gamma(E)$, by
${\cal S}(E;h)=S_0(E)+hS_1(E)+{\cal O}(h^2)$ with
$$\leqalignno{
&S_0(E)=\int_{\gamma_E}\xi\,dx&(1.6)\cr
&S_1(E)=-\int_0^{T(E)} H_1(y(t),\eta(t))\,dt+{1\over2i}\sum_{j=1}^d\mu_j(E)+\pi{g_\ell\over2}&(1.7)\cr
}$$
Here $\mu_j(E)=\mu_j+{\cal O}(E)$ is Floquet exponent at energy $E$, 
$g_\ell\in{\bf Z}$ Gelfand-Lidskiy or Cohnley-Zehnder index of $\gamma(E)$ (depending
only on elliptic elements). Then 
the resonances of $H$ associated with the family $\gamma(E)$ for $E$ in ${\cal W}_h=[E_0-\e _0,E_0+\e _0]-i]0,h\log(1/h)]$
are given by the generalized BS quantization condition 
$${1\over2\pi h}{\cal S}(E;h)+{1\over2i\pi}\bigl(\sum_{j=1}^d hk_j\mu_j(E)+{\cal O}(h^2|k|^2)\bigr)=mh, \ m\in{\bf Z}, \ k\in{\bf N}^d\leqno(1.8)$$
provided $|m|h\leq \e_0$, $|k|\leq \const \log(1/h)$. 
\smallskip
This remains true, at the price of technical difficulties, when replacing the the width $h\log(1/h)$ of $W_h$ 
by $h^{1-\delta}$, $0<\delta<1$. We stress that this theorem says in general nothing about other resonances described in Theorem 1.1,
unless $\gamma_0$ is purely hyperbolic, in which case the periodic orbits $\gamma(E)$ are isolated, and thus we can assume 
${\cal K}(E)=\gamma(E)$.  
\medskip
\noindent{\it e) Remarks on the trace formulas}
\smallskip
In the self-adjoint case (e.g. the geodesic flow on a compact manifold with negative curvature) 
trace formulas have been considered for hyperbolic or semi-hyperbolic flows. They are expressed 
in the time variable $t$
(trace of the propagator or wave group, see [Zel]), or
in the energy
variable $E$ (trace of the semi-classical Green function, see [Vo] and references therein). 

In case $H$ is the geodesic flow on a compact Riemannian manifold $(M,g)$, Zelditch [Zel2] computed the singular part
of the trace of the wave group $U(t)$.  It is obtained as a term $e_0(t)$ (involving the fixed points of the flow), plus
the sum over all periodic geodesics $\gamma$ on $M$, of ``wave trace invariants''
$e_\gamma(t)$, using non commutative residues, 
that can be computed as an asymptotic series (asymptotics with respect to smoothness). This formula does not
involve other periodic orbits (clustering on each $\gamma$ when $\gamma$ is of semi-hyperbolic type), but their contribution
would appear when investigating ``convergence'' (in the sense of resurgence) of the series defining $e_\gamma(t)$. 
Note that the semi-classical parameter is obtained in scaling the variables microlocally near a periodic geodesic
to bring the Hamiltonian in BNF. 
The same situation
is likely to appear for resonances in the non compact case. 

The trace of the semi-classical Green function instead, near some fixed $E$  can be expressed formally by a sum of terms $R_\gamma(E)$ 
labelled by the classical periodic orbits $\gamma$ having energy $E$. The poles of $R_\gamma(E)$ are 
precisely localized by an implicit equation such as Bohr-Sommerfeld quantization condition of Theorem 1.2:
since it would give complex energies, it is called by Voros the ``generalized quantization condition'', to stress that 
periodic orbits are not necessarily associated with bound states. 
This paradox could be settled in the framework or resurgence theory. 

In the context of resonances the paradox of complex poles disappears. Thus it 
would be tempting to look for a trace formula as in [SjZw].
First we recall Helffer-Sj\"ostrand formula [DiSj]. Let $H$ be a self-adjoint operator, 
$\chi\in C^2_0({\bf R})$ and $\widetilde \chi\in C^1_0({\bf C})$ an almost analytic extension of $\chi$ satisfying
$\overline\partial\widetilde \chi(z)={\cal O}(|\im z|)$. Then
$$\chi(H)=-{1\over\pi}\int\overline\partial\widetilde \chi(z)(z-H)^{-1}\,L(dz)\leqno(1.9)$$
Following the remark after the proof of Theorem 8.1 in [DiSj], this could be generalized to
classes of non selfadjoint operators, and applied to $H_\theta$. Then we may
use (1.9) as a definition. 

So let $f_N\in C^\infty({\bf R})$ be such that $\supp \widehat f_N\subset]0, T_0 N[$, $\phi_N\in C_0^\infty({\bf C})$ whose support in $\im z$
has to be chosen suitably near 0 (of width ${\cal O}(h)$), $A_N(y,hD_y)$ be a $h$-PDO cutoff equal to 1 in a small neighborhod of $\gamma(z)$, 
$M(z)$ be the monodromy operator computed either in Sect.2 or Sect.3, and $k(N)\to\infty$. In the case of resonances,
it is plausible to expect a ``trace formula''
modelled after this of [SjZw], namely  
$$\eqalign{
&\Tr f_N(H_\theta/h) \phi_N(H_\theta;h) A_N(y,hD_y)\approx\cr
&{1\over2i\pi}\Sum_{k=1}^N\Tr \int_{\bf R}f_N(z/h)\phi_N(z;h) M(z)^{k-1}{dM\over dz}(z)\,dz+{\cal O}(h^{k(N)})
}$$
just keeping positive values of $N$ to account for the time reversal symmetry breaking.
This however, seems again far from 
reach, especially because resonances proliferate near the real axis in $W_h$ as $h\to0$. 

At last we note that in the framework of resonance scattering outside convex obstacles, trace formulas (or the related zeta function) 
in the energy representation 
are given by Ikawa [Ik], and the situation is better understood. It is similar to our case, when Poincar\'e map
has no elliptic element. 
\medskip
\noindent{\it Acknowlegments}: We are grateful to Alain Chenciner, Sergey Bolotin and Andr\'e Voros for useful information, 
and to a referee for useful remarks. 
The second author was partially supported by Grant 
PRC No. 1556 CNRS-RFBR 2017-2019 
``Multi-dimensional semi-classical problems of 
Condensed Matter Physics and Quantum Mechanics''.
\bigskip
\noindent {\bf 2. A hint on the proof of Theorem 1.1}.
\smallskip
\noindent {\it a) The (absolute) monodromy operator}
\smallskip
The energy parameter will be denoted by $z$, and for the moment we work at a formal level, i.e. $H(y,hD_y;h)$ denotes the self-adjoint operator. 
We shall follow mainly [SjZw], making use of 2 equivalent Poincar\'e sections,
but taking care eventually of the (semi-) hyperbolic structure of the flow.
   
Before entering the actual constructions, we recall how to define the monodromy operator and solve Grushin problem
in the simple situation (see [SjZw], [IfaLouRo]), 
where $H=hD_t$ acts on $L^2({\bf S}^1)$ with periodic boundary condition
$u(t)=u(t+2\pi)$. 
Solving for $(H-z)u(t)=0$, we get two solutions with the same expression but defined on different charts
$$u^{a}(t)=e^{izt/h}, -\pi<t<\pi, \quad u^{a'}(t)=e^{izt/h}, 0<t<2\pi\leqno(2.13)$$
indexed by angles $a=0$ and $a'=\pi$ on ${\bf S}^1$. All angles will be computed mod $2\pi$. 
In the following we take advantage of the fact that 
these functions differ but when $z$ belongs to the spectrum of $H$. 

Let also $\chi^a\in C_0^\infty({\bf S}^1)$
be equal to 1 near $a$, $\chi^{a'}=1-\chi^a$. To fix the ideas, we may assume that $\chi^a$ drops down to 0 near $t=-{\pi\over2}$ and $t={\pi\over2}$.
(which belongs to both charts). We set $F^a_\pm={i\over h}[P,\chi^a]_\pm u^a$, where $\pm$ denotes the part of the 
commutator supported in the half circles $0<t<\pi$ and $-\pi<t<0$ mod $2\pi$. Similarly  
$F^{a'}_\pm={i\over h}[P,\chi^{a'}]_\pm u^{a'}$, and we may assume that $\chi^{a'}$ drops to 0 near $t={\pi\over2}$ and $t={3\pi\over2}$ mod $2\pi$.
Modulo ${\cal O}(h^\infty)$ (as all constructions in this work, so we shall not dwell on this anymore)
distributions $F^{a}_\pm,F^{a'}_\pm$ do not depend on the choice of $\chi^a$ above
since we may expand the commutator when applying to distributions defined on a single chart (2.13) and use that $H$ is self-adjoint. 
\medskip
\noindent {\it Remark}: It is convenient to view $F^a_+-F^a_-$ and $F^{a'}_+-F^{a'}_-$ as belonging to co-kernel of $H-z$ in the sense they are not
annihilated by $H-z$. If we form Gram matrix 
$$G^{(a,a')}(z)=\pmatrix{(u^a|F^{a}_+-F^{a}_-)&(u^{a'}|F^{a}_+-F^{a}_-)\cr (u^a|F^{a'}_+-F^{a'}_-)&(u^{a'}|F^{a'}_+-F^{a'}_-)\cr}$$
an elementary computation shows that 
$\det G^{(a,a')}(z)=-4\sin^2(\pi z/h)$,
so the condition that $u^a$ coincides with $u^{a'}$ is precisely that $z=kh$, with $k\in{\bf Z}$ (see [IfaRLouRo] for details,
where a convention slightly different from [SjZw] has been made).
\smallskip

Starting from the point $a=0$ we associate with $u^a$ the multiplication operator $v_+\mapsto I^a(z)v_+=u^a(t)v_+$ on ${\bf C}$,
i.e. Poisson operator with  
``Cauchy data'' $u(0)=v_+\in{\bf C}$. 
Similarly multiplication by $u^{a'}$ defines Poisson operator $I^{a'}(z)v_+=u^{a'}(t)v_+$, 
which another ``Cauchy data'' $v_+$ at $a'=\pi$.
\smallskip
Now we turn to the general case. The situation of [NoSjZw] simplifies 
since we need only 2 equivalent Poincar\'e sections (modulo moving around $\gamma_0$).
Let $m_0\in\gamma_0$, and $m_1=\exp {T_0\over2}X_{H}(m_0)$ be two ``distinguished'' points on $\gamma_0$,
which play respectively the role of $t=\pi/2$, $t=-\pi/2$ mod $2\pi$ in the example above. 
Assume for simplicity that $H(y,hD_y;h)$ verifies time-reversal symmetry, so that $m_0$ is again the point along $\gamma_0$
reached from $m_1$ within time $T_0/2$, where $T_0$ is the (primitive) period of $\gamma_0$. The corresponding cut-off
near the orbit will still be denoted by $\chi^a$ and $\chi^{a'}$. 

Since the energy shell $E=0$ is non critical, near every $m\in\gamma_0$, $H(y,hD_y;h)$ 
can be reduced microlocally to $hD_{t}$, i.e. there exists a local canonical transformation $\kappa: T^*M\to {\cal U}\subset T^*{\bf R}^n$
defined near $((0,0),m)$, and a $h$-FIO ${\cal T}$, associated with the graph of $\kappa$, elliptic near $((0,0),m)$, such that
${\cal T}H(y,hD_y;h)=hD_{t}{\cal T}$ near $((0,0),m)$.

Fix $m$, and construct a corresponding ${\cal T}$; 
if we define 
$$\ker _{m}(H)={\cal T}^{-1}\ker(hD_{t})$$ 
we can identify $\ker _{m}(H)$ with semi-classical distributions on ${\bf R}^d$
(i.e. on a Poincare section) microlocally near $(m,(0,0))$; we denote this identification by $K:{\cal D}'({\bf R}^d)\to \ker _{m}(H)$.  
Now we solve $(hD_{t}-z)u(x)=0, u|_{t=0}=v(x')$ in ${\cal D}'({\bf R}^n)$ by $u(x)=e^{izt/h}v(x')$, $x=(t,x')$. As in the Example,
we obtain this way two Poisson operators $v_+\mapsto I^a(z)v_+=e^{izt/h}v_+$ when $-\pi<t<\pi$ (forward) and 
$v_+\mapsto I^{a'}(z)v_+=e^{izt/h}v_+$ when $0<t<2\pi$ (backward), defined on ``Cauchy data'' $v_+\in{\cal D}'({\bf R}^d)$.
Working locally, we can ignore their domain, 
and call them both $K(z)$, but moving along the flow in either direction, we introduce new canonical charts $(\kappa,{\cal U})$
and construct new FIO's ${\cal T}$ accordingly. By compactness, we can cover $\gamma_0$ with a finite set of such $(\kappa,{\cal U})$. 
Assume the intersection of Poincar\'e sections with the domain of definitions of ${\cal T}$ contains (strictly) the trapped set. 

Instead of specifying $\gamma_0$, 
it is more convenient to select the orbit $\gamma(z)$, which is periodic with respect to the Hamilton flow of $H$ at energy $z$,
with period $T(z)$. Accordingly, we change $m_0,m_1$ to $m_0(z),m_1(z)$. All $\gamma(z)$ are mapped diffeomorphically to ${\bf S}^1$
by the Hamilton flow, so moving once around $\gamma(z)$ means moving once around ${\bf S}^1$ in the Example.

Varying $m$ on the orbit $\gamma(z)$, we obtain the forward/backward extensions (standing for $u^a,u^{a'}$ in (2.13)), 
independent of $m(z)\in\gamma(z)$
$$I_\pm(z):{\ker }^\pm _{m(z)}(H-z)\to{\ker }_{\gamma(z)}(H-z)\leqno(2.18)$$ 
where ${\ker }^\pm _{m(z)}(H-z)$ denotes the space of forward/backward solutions near $m(z)$.
Operators $I_\pm(z)$ 
are (microlocally) injective. 
Thus we obtain the exact sequence (with obvious notations)
$$0\longrightarrow{\ker }^+ _{m(z)}(H-z)\oplus{\ker }^- _{m(z)}(H-z)\longrightarrow
{\ker }_{\gamma(z)}(H-z)\longrightarrow0\leqno(2.19)$$
where the 2nd arrow is $I_+(z)\oplus I_-(z)$ and the 3rd arrow is $H-z$. (2.19) remains true if we change 
${\ker }^+ _{m(z)}(H-z)\oplus{\ker }^- _{m(z)}(H-z)$, e.g. to ${\ker }^+ _{m_0(z)}(H-z)\oplus{\ker }^- _{m_1(z)}(H-z)$. 
Let $K_f(z)=I_+(z)K(z)$,  $K_b(z)=I_+(z)K(z)$. Since $K(z)$ identifies microlocally ${\ker }^\pm _{m(z)}(H-z)$
with ${{\cal D}'}_{m(z)}({\bf R}^d)$ (semi-classical distributions microlocalized on Poincar\'e section at $m(z)$), we have also
$$0\longrightarrow{{\cal D}'}_{m(z)}({\bf R}^d)\oplus{{\cal D}'}_{m(z)}({\bf R}^d)\longrightarrow
{\cal D}'_{\gamma(z)}({\bf R}^n)\longrightarrow0\leqno(2.20)$$
where the 2nd arrow is $K_f(z)\oplus K_b(z)$.
Let $W_+,W_-$ be nghbhds of $m_0,m_1$ in $T^*M$, such that 
$$\{m:\chi^a(m)\neq1\}\cap\{m:\chi^a(m)\neq0\}=\{m:\chi^{a'}(m)\neq1\}\cap\{m:\chi^{a'}(m)\neq0\}\subset W_-\cup W_+\leqno(2.21)$$
Following [SjZw], we define microlocally near $W_+$ the 
{\it (absolute) quantum monodromy operator}
${\cal M}(z):{\ker }^+_{m_0(z)}(H-z)\to {\ker }^-_{m_0(z)}(H-z)$ by
$$I_+(z)f=I_-(z){\cal M}(z)f, \ f\in \ker _{m_0(z)}(H-z), \ \hbox{microlocally near} \ W_-\leqno(2.23)$$
Clearly, we can interchange the roles of $m_0(z)$ and $m_1(z)$ in this definition: ${\cal M}(z)$ is independent of the section. 
We define the {\it quantum monodromy operator} $M(z):{\cal D}'({\bf R}^d)\to{\cal D}'({\bf R}^d)$ as follows.
Let
$$\widetilde B(z)=K(z)^*{i\over h}[H,\chi^a ]_{W_+} K(z):{\cal D}'_{m_0}({\bf R}^d)\to{\cal D}'_{m_0}({\bf R}^d)$$
Following [SjZw], we check that $\widetilde B(z)$ is a $h$-PDO, defined  microlocally near $(0,0)\in T^*{\bf R}^d$,
positive and formally self-adjoint. So
$L(z)=K(z)\widetilde B(z)^{-1/2}$ verifies the ``flux norm'' identity
$$L(z)^*{i\over h}[H,\chi^a ]_{W_+} L(z)=\Id \ \hbox{microlocally near} \ (0,0)\in T^*{\bf R}^d\leqno(2.25)$$
Using (2.20), we see that $K(z):{{\cal D}'}_{m_0(z)}({\bf R}^d)\to \ker _{m_0(z)}(H-z)$ is invertible, so is $L(z)$ and its inverse is 
$$L(z)^{-1}=R_+(z)=L(z)^*{i\over h}[H,\chi^a ]_{W_+}\leqno(2.26)$$
We call the {\it quantum monodromy operator} 
$$M(z)=L(z)^{-1}{\cal M}(z)L(z)\leqno(2.28)$$
This is a $h$-FIO, whose canonical relation is precisely the graph of Poincare map ${\cal P}_0(z)$, i.e. 
$$M(z)\in I^0\bigl({\bf R}^d\times{\bf R}^d;C'), \ C'=\{(x,\xi,x',-\xi'): (x,\xi)={\cal P}_0(z)(x',\xi')\}$$
We make more precise [SjZw,Prop.4.5] (still before any analytic dilation), taking also into account the hyperbolicity of the flow.
Recall the flow is expanding in some direction of Poincar\'e section, and contracting in the orthogonal one (for the symplectic structure).
Denote by $D(W_+)\subset W_+$ ($D$ like {\it departure}) a neighborhood of the outgoing manifold in $W_+$ and 
$A(W_+)\subset W_+$ ($A$ like {\it arrival}) a neighborhood of the incoming  manifold in $W_+$ (see [NoSjZw], [NoZw]). 
By the same letter we denote the space
of distributions microlocalized near that set. 
\medskip
\noindent {\bf Proposition 2.1}: For real $z$, the monodromy operator $M(z)$ is microlocally ``unitary''
$D(W_+)\to A(W_+)$, and similarly
for complex $z$, in the sense that the adjoint of $M(z)$ is equal to $\bigl(M(\overline z)\bigr)^{-1}$. 
\smallskip
\noindent {\it Proof}: Let $v\in{\cal D}'_{m(z)}({\bf R}^d)$ microlocally supported near (0,0),
and $u=L(z)v\in {\cal D}'_{\gamma(z)}({\bf R}^n)$, we compute (dropping the variable $z$ from the notations)
$(Mv|Mv)=(L^{-1}{\cal M}u|L^{-1}{\cal M}u)$. 
By inserting (2.26) on the left of the scalar product we get
$$(Mv|Mv)=\bigl({i\over h}[H,\chi^a ]_{W_+}{\cal M}u|{\cal M}u\bigr)_{L^2({\bf R}^n)}=
\bigl((I_-)^{-1*}{i\over h}[H,\chi^a ]_{W_+}{\cal M}u|I_-{\cal M}(z)u\bigr)$$
where we have also introduced the backward extension operator $I_-=I_-(z)$ as in (2.18). Next we have, for $\delta>0$ sufficiently small, 
and $0<t<T(0)/2+\delta$,
${i\over h}[H,\chi]_{W_+}=I_-^*{i\over h}[H,\chi^t ]_{W_+^t}I_-$, where $\chi^t =\chi\circ\exp -tX_H\equiv \chi^a(\cdot -t)$, 
and $W_+^t=\exp(-tX_H)W_+$, which corresponds to moving $\supp \chi$ in the direction opposite to the flow of $X_H$, and 
$W_+$ simultaneously so that (2.21) holds. Hence
$$(Mv|Mv)=\bigl({i\over h}[H,\chi^t ]_{W_+^t}I_-{\cal M}u|I_-{\cal M}u\bigr)_{L^2({\bf R}^n)}\leqno(2.31)$$
Similarly, inserting (2.25) on the left of the scalar product $(v|v)=(L^{-1}u|L^{-1}u)$ we get 
$$(v|v)=\bigl({i\over h}[H,\chi^{-t} ]_{W_+^{-t}}I_+u|I_+u\bigr)\leqno(2.32)$$
For $t\sim T(0)/2$, and $\WF _hv$ sufficiently close to (0,0) so that $\WF _h{i\over h}[H,\chi^{-t} ]_{W_+^{-t}}I_+u\subset W_-$, we get
$(v|v)=\bigl({i\over h}[H,\chi^{t} ]_{W_+^{t}}I_-{\cal M}u|I_-{\cal M}u\bigr)$, and comparing (2.31) with (2.32) gives $(Mv|Mv)=(v|v)$.
The Proposition follows easily from the definition of $D(W_+),A(W_+)$. $\clubsuit$ 
\smallskip
To fully restore ``unitarity'' of $M$, so that Grushin problem be well-posed, we need to introduce the weighted
Sobolev spaces, or/and the complex Lagrangian deformations. Let us conclude by writing $M(z)$ in a form similar to [NoSjZw,(4.33)].
This is done in 2 steps: let $K_0(z)=K(z),K_1(z)$ be Poisson operators at $m_0(z),m_1(z)$, and $L_0(z),L_1(z)$ be
the normalized ones. The monodromy operator
from $m_0(s)$ to $m_1(z)$ is $M_{0,1}(z)=L_1(z)^*{i\over h}[H,\chi^a ]_{W_+} L_0(z)$, and this 
from $m_1$ back to $m_0(z)$, $M_{1,0}(z)=L_0(z)^*{i\over h}[H,\chi^a ]_{W_-} L_1(z)$. Then 
$$M(z)=M_{1,0}(z)M_{0,1}(z)\leqno(2.33)$$
This will simplify (for the simplified problem) in action-angle coordinates as we shall see later.
\medskip
\noindent {\it b) Intertwining ${\cal M}(z)$ with ${\cal M}(w)$}.
\smallskip
Structural stability for hyperbolic flows ([KaHa,Thm.18.2.3]) recalled in Sect.1 carries to the monodromy operator. 
Namely, following [SjZw], let $p(z)=H_0-z$. We call a {\it classical time function} a solution $q(z)$ (which can be chosen independent of $z$) of 
$$X_{p(z)}q(z)=\{p(z),q(z)\}=1$$
(Lie differentiation). Thus $(q(z),p(z))$ are just the restriction to $\gamma(z)$ (in the energy shell $p(z)=0$)
of (symplectic) Darboux coordinates $(t,z)$ along $\gamma(z)$, 
adapted to the Stable/Un- stable/Center manifold.
 Since $q(z)$ is a multi-valued function, we call {\it first return classical time function}, and denote by $q_\partial(z)$
its continuation to the second sheet. Thus we have, with a slight abuse 
of notations
$$q_\partial(z)\bigl(m(z)\bigr)=q\circ\exp T(z)X_{p(z)}\bigl(m(z)\bigr)$$
and
$$\bigl(q_\partial(z)-q(z)\bigr)|_{\gamma(z)}=T(z)\leqno(2.38)$$ 
where $T(z)={dJ\over dz}$, $J(z)=\int_{\gamma_z}\eta\,dy$ being the classical 
action along $\gamma(z)$.   
We call a {\it quantum time} (resp. {\it first return quantum time} 
a solution $Q(z)$, in the $h$-PDO's sense, of
$$\Id={i\over h}[P(z),Q(z)], \quad \Id={i\over h}[P(z),Q_\partial(z)]\leqno(2.39)$$
with principal symbols $q(z),q_\partial(z)$ respectively. Here $P(z)=H-z$. In the case $P(z)$ is self-adjoint, we can assume $Q(z)$ 
and $Q_\partial(z)$ are self-adjoint (here again we work formally, but we shall need to take hyperbolicity into account as before).  We have
$$Q_\partial(z)-Q(z):\ker _{m(z)}(H-z)\to \ker _{m(z)}(H-z)$$
Next we construct $h$-FIO's that will intertwine Poisson operators 
at different energies, and consider the following system of equations
$$\eqalign{
&\bigl(hD_z-Q(z)_L\bigr)U(z,w)=hD_zU(z,w)-Q(z)U(z,w)=0\cr
&\bigl(hD_w+Q(w)_R\bigr)U(z,w)=hD_zU(z,w)+U(z,w)Q(w)=0\cr
}\leqno(2.41)$$
with initial condition $U(0,0)=\Id$. We can write (2.41) as ${\cal L}_{L/R}(z,w)U(z,w)=0$, 
with ${\cal L}_L=hD_z-Q(z)_L$, ${\cal L}_R=hD_w+Q(w)_R$, and the solvability condition is ensured by
the commutation relation $[{\cal L}_L,{\cal L}_R]=0$. 
It turns out that $U(z,w)$ can be constructed in the class of $h$-FIO's on ${\bf R}^n$, microlocally near $m\in\gamma$. 
For the model, $U(z,w)$ is just the multiplication operator by $e^{it(z-w)/h}$. 
We notice that (2.41) implies $U(z,z)=\id$, $U(z,w)U(w,v)=U(z,v)$, and $U(w,z)^*=U(z,w)$ when $H(y,hD_y;h)$ is self-adjoint
We have $K(z)=U(z,w)K(w)$, and differentiating gives $hD_zK(z)=Q(z)K(z)$.
Further, varying $m$, we extend $U(z,w)$ in the forward and backward regions, to $U_\pm(z,w)$. We have
$$(H-z)U_\pm(z,w)=U_\pm(z,w)(H-w), \quad U(z,w)I_\pm(z)=I_\pm(w)\leqno(2.42)$$
Changing $Q(z)$ to $Q_\partial(z)$ in (2.41), we can solve for $U_\partial(z,w)$ with same properties as $U(z,w)$. There follows the 
\medskip
\noindent{\bf Proposition 2.2}: We have the intertwining property
$${\cal M}(z)U(z,w)=U_\partial(z,w){\cal M}(w)$$
and the quantum monodromy operator satisfies the equation
$$hD_zM(z)=K(z)^{-1}(Q_\partial(z)-Q(z))K(z)M(z)$$
\medskip
\noindent {\it c)  Grushin problem}
\smallskip
Consider again the model case, with the notations of Sect.2. Introduce the ``trace operator''
$R_+(z)u=u(0)$, if $u(t)=e^{izt/h}v$ with $v_+=u(0)$, we check that 
$$K(z)^*{i\over h}[P,\chi^a]u=\int_{-\pi}^\pi e^{-izt/h}(\chi^a)'(t)u(t)\,dt=v_+$$
Consider also the multiplication operators
$$E_+(z)=\chi^a I^{a}(z)+(1-\chi^a)I^{a'}(z), \ R_-(z)={i\over h}[P,\chi^a]_-I^{a'}(z), \ E_{-+}(z)=1-e^{2i\pi z/h}$$
We claim that 
$${i\over h}(P-z)E_+(z)+R_-(z)E_{-+}(z)=0\leqno(2.51)$$
Namely, evaluating on $0<t<\pi$, 
we have $I^a(z)=e^{itz/h}, I^{a'}(z)=e^{itz/h}$, while evaluating on $-\pi<t<0$, 
$I^a(z)=e^{itz/h}, I^{a'}(z)=e^{i(t+2\pi)z/h}$. Now 
${i\over h}(P-z)E_+(z)=[P,\chi^a]\bigl(I^a(z)-e^{i\pi z/h}I^{a'}(z)\bigr)$ vanishes on $0<t<\pi$, while is equal to
$R_-(z)E_{-+}(z)$ on $-\pi<t<0$. So (2.51) follows. 
Hence Grushin problem
$${\cal P}(z;h){u\choose u_-}=\pmatrix{{i\over h}(P-z)&R_-(z)\cr R_+(z)&0}{u\choose u_-}={v\choose v_+}\leqno(2.52)$$
with $v=0$ has a solution $u=E_+(z)v_+$, $u_-=E_{-+}(z)v_+$, and $E_{-+}(z)$ is the effective Hamiltonian. As we show below,
we can find $E(z)$ such that problem (2.52) is well posed, ${\cal P}(z)$ is invertible, 
and 
$${\cal P}(z)^{-1}=\pmatrix{E(z)&E_+(z)\cr E_-(z)&E_{-+}(z)}\leqno(2.54)$$ 
with
$$(P-z)^{-1}=E(z)-E_+(z)E_{-+}(z)^{-1}E_-(z)\leqno(2.55)$$
In our case however, because of hyperbolicity, we need to introduce the weighted spaces (or Lagrangian deformations)
so that (2.52) be well-posed. Still we start to proceed within the formalism of Sect.2. 
Recall $R_+(z)$ from (2.26). So if $v\in{\cal D}'({\bf R}^d)$,
$u=L(z)v$ solves near any $m(z)$ 
$$(H-z)u=0, \quad R_+(z)u=v\leqno(2.57)$$
To obtain a Cauchy problem globally near $\gamma(z)$, we need to introduce $R_-(z)$. Recall $K_{f/b}(z)=I_\pm(z)K(z)$, 
which we normalize to $L_{f/b}(z)=I_\pm(z)L(z)$ as in (2.26). By (2.23) and (2.28), we have
$$L_f(z)=L_b(z)M(z), \ \hbox{microlocally near } \ W_-\times (0,0)\leqno(2.58)$$
and solve (2.57) in $\Omega\setminus W_-$ ($\Omega$ neighborhood of $\gamma(z)$) as in the argument after (2.51) by
$$E_+(z)v_+=\chi^a L_f(z)v_++(1-\chi^a)L_b(z)v_+$$
so that in particular $E_+(z)v_+=L(z)v_+$ in $W_+$ (since $L_f(z)=L_b(z)$ in $W_+$), and 
$$R_+(z)E_+(z)v_+=L(z)^*{i\over h}[H,\chi^a]_{W_+}L(z)v_+=v_+$$
by (2.25). Applying $H-z$, using (2.58) and $(H-z)E_+(z)v_+=0$ in $W_+$, we find that,
with $R_-(z)={i\over h}[H,\chi^a]_{W_-}L_b(z)$, and $E_{-+}(z)=\id -M(z)$, $u=E_+(z)v_+$, $u_-=E_{-+}(z)v_+$ solve (formally)
the problem ${\cal P}(z){u\choose u_-}={0\choose v_+}$ near $\gamma(z)$. This implies that the microlocal inverse of ${\cal P}(z)$
should be of the form 
${\cal E}(z)=\pmatrix{E(z)&E_+(z)\cr E_-(z)&E_{-+}(z)}$, and we still have to find $E(z), E_-(z)$. So we try to solve the inhomogeneous problem
${\cal P}(z){u\choose u_-}={v\choose v_+}$ near $\gamma(z)$, and introduce the forward/backward fundamental solutions of $H-z$, namely
$E_f(z)=\int_0^\infty e^{-it(H-z)/h}\,dt$, $E_b(z)=\int_{-\infty}^0 e^{-it(H-z)/h}\,dt$, which of course assume a simple form after taking $H$
microlocally to $hD_t$. The construction of $E(z)$ is more involved (see [SjZw], [NoSjZw]), but an argument like in Proposition 2.1 leads to
$$E_-(z)=-\bigl(M(z)K_f(z)\chi+K_b(z)^*(1-\chi)\bigr)$$
Next we need to specify the right spaces where Grushin problem is well posed. This is done by introducing microlocal weights as in the Appendix,
encoding the trapped set. We eventually get Theorem 1.1 as in [NoSjZw]; details will be given elsewhere.
\medskip
\noindent {\bf 3. An ``approximate'' theory}.
\smallskip
Here we ``neglect'' the occurrence of infinitely many periodic orbits near $\gamma_0$. 
It is plausible that this theory would still provide a good description of the resonant spectrum close to the real axis,
since orbits with large period are quite unstable and contribute to the spectrum only far away from the real axis.
Moreover, it becomes exact in the particular case where there are no elliptic elements,
because such periodic orbits are isolated.  
At last, it provides BS quantization rules for the family $\gamma(E)$, which are known to hold also in the semi-hyperbolic case.

Using complex coordinates, we may also reduce the center manifold $C$ to $\overline\gamma$ by moving the elliptic subspaces into $N_\pm$.

\medskip
\noindent {\it a) Birkhoff normal form}
\smallskip
Our approach relies on the classical BNF for the principal symbol $H_0$ of $H$. The first step takes $H_0$ to the form
$H_0(y,\eta)=-\tau+\langle B_0x,\xi\rangle+g(\tau)+{\cal O}(|\tau,|x,\xi|^2|^2)$
(the natural orientation of $\gamma_0$ has been reversed). Here $(t,\tau)$ parametrize $T^*\overline\gamma$, $(x,\xi)$ are transverse
variables on Poincare section, 
$g(\tau)=\tau+f(-\tau)={\cal O}(\tau^2)$, and $f$ parametrizes energy
according to $f(-\tau)=E$; it is related to the period $T(E)$ of $\gamma(E)$ by $f'(-\tau)={2\pi\over T\circ f(-\tau)}$, with $f'(0)=1$. 
\medskip
\noindent {\bf Proposition 3.1} [Br],[GuPa]: Assume that Floquet exponents satisfy the strong non-reson-ance condition (H.6). 
Then in a nghbhd of $\gamma_0$, there exists symplectic coordinates $(t,\tau,x,\xi)$, $t\in[0,2\pi]$, such that for all $N\geq1$, 
we can find a canonical transformation $\kappa_N$ with 
$$H_0\circ\kappa_N=-\tau+\Sum_{j=1}^d\mu_jQ_j(x,\xi)+H_0^{(N)}(\tau;Q_1,\cdots,Q_n)+{\cal O}\bigl(|\tau,|x,\xi|^2|^{N+1}\bigr)\leqno(3.1)$$
where $H_0^{(N)}(\tau;Q_1,\cdots,Q_d)={\cal O}(|\tau,Q|^2)$ is a polynomial of degree $N$, and the remainder term
${\cal O}\bigl(|\tau,|x,\xi|^2|^{N+1}\bigr)$ is $2\pi$-periodic in $t$.
Here $\mu_jQ_j(x,\xi)$ is a polynomial of the form ${\mu_j\over2}(\xi_j^2-x_j^2)$ (or $\mu_jx_j\xi_j$) (hr element), 
$-{i\mu_j\over2}(\xi_j^2+x_j^2)$
(ee element),
$c_j(x_{2j-1}\xi_{2j-1}+x_{2j}\xi_{2j})-d_j(x_{2j-1}\xi_{2j}-x_{2j}\xi_{2j-1})$, $\mu_j=c_j+id_j$
(hc or loxodromic elements) which also take the form $\mu_jx_j\xi_j$ in complex coordinates. 
\smallskip
This BNF carries to the semi-classical setting (see also [Zel] for high energy expansions):
\medskip
\noindent {\bf Proposition 3.2} [GuPa]: Under hypotheses above, conjugating with a $h$-FIO microlocally unitary near $\gamma_0$,
$H^w(y,hD_y;h)$ can be taken formally to 
$$\eqalign{
H^{(N)}&(hD_t,x,hD_x;h)=-hD_t+\Sum_{j=1}^d\mu_jQ_j^w+H_0^{(N_0)}\bigl(hD_t;Q_1^w,\cdots,Q_d^w\bigr)+\cr
&+hH_1^{(N_1)}\bigl(hD_t;Q_1^w,\cdots,Q_d^w\bigr)+\cdots\cr}
\leqno(3.2)$$
as a polynomial depending on the $n$ ``variables'' $(hD_t,Q_j^w)$, with for instance when $\mu_j$ is real,
$$Q_j^w={1\over2}(x_jhD_{x_j}+hD_{x_j}x_j)=\Op ^wQ_j$$ 
where the \dots stand for terms ${\cal O}(h^2)$, as well as operators with coefficients ${\cal O}(h^\infty)$ and periodic in time.
$N_j$ denotes the order of expansion as a Birkhoff series of Hamiltonian $H_j$, and $N=(N_0,N_1,\cdots)$ any sequence of integers.   
Moreover, allowing for complex coordinates, one can formally assume that $Q_j^w={1\over2}(x_jhD_{x_j}+hD_{x_j}x_j)$ for all types of elements
(ee or he). 
\smallskip
Keeping the leading part in (3.2) the Model Hamiltonian,  
$$H_{\mod }(hD_t,x,hD_x;h)=-hD_t+\Sum_{j=1}^d\mu_jQ_j^w(x,hD_x)$$ 
with periodic boundary conditions on ${\bf S}^1\times{\bf R}^d$
serves as a guide-line as $hD_t$ did in Sect.2. 
\medskip
\noindent {\it b) Microlocalisation in the complex domain}
\smallskip
Taking into account that there exists an escape function
outside the trapped set $\gamma(E)$, the 
most relevant region of phase-space for such deformations is a neighborhood of $\gamma(E)$.
Here we make a complex scaling of the form $(x,\xi)\mapsto(e^{i\theta}x,e^{-i\theta}\xi)$ (independent of $E$), followed also by a
deformation in the $(t,\tau)$ variables. Rather then using weighted spaces as in Sect.2, our main tool is the method of Lagrangian
deformations. Namely we perform a
FBI transformation (metaplectic FIO with complex phase) which takes the form, in coordinates 
$(s,y;t,x)\in T^*{\bf R}^n\times T^*{\bf C}^n$  adapted to $\Gamma_\pm$ as in BNF
$$T_0u(x,h)=\int e^{i\varphi_0(t,s;x,y)/h}u(s,y)\,ds\,dy, \ u\in L^2({\bf R}^n)$$ 
where $\varphi_0(t,s;x,y)=\varphi_1(t,s)+\varphi_2(x,y)$,
$\varphi_1(t,s)={i\over2}(t-s)^2$, $\varphi_2(x,y)={i\over2}\bigl[(x-y)^2-{1\over2}x^2\bigr]$.
The corresponding pluri-subharmonic (pl.s.h.) weight is $\Phi_0=\Phi_1+\Phi_2=(\im t)^2/2+|x|^2/4$.
In a very small neighborhood of $\gamma(E)$, whose size
will eventually depend on $h$, corresponding to $\theta=-\pi/4$, and that we call the ``phase of inflation'',
$T_0H^w(y,hD_y;h)T_0^{-1}=\widetilde H(hD_t,x,hD_x;E,h)$ assumes BNF and is approximated at leading order by the Model Hamiltonien. 
In a somewhat larger neighborhood of $\gamma(E)$, which we call
the ``linear phase'', we choose $\theta$ small enough, and get a new pl.s.h. weight $\widetilde\Phi_\theta(t,x)$.
Farther away from $\gamma(E)$ (in the ``geometric phase'') the weight is implied by the 
escape function. All these weights are
patched together in overlapping regions,
so to define a globally pl.s.h. function in complex $(t,x)$ (or $y$) space. It
determines the contour integral for writing realizations of $h$-FIO's in the complex domain [Sj]
in $H_\Phi$ spaces, conjugating $H^w(y,hD_y;h)$ to 
a $h$-PDO everywhere elliptic but on $\overline\gamma$. In particular near $\overline\gamma$
$$|(H-E)|_{\Lambda_{\widetilde\Phi}}\sim|x|^2+|\im t-\tau|\leqno(3.3)$$
\medskip
\noindent {\it c) Poisson operator, its normalisation and the monodromy operator}
\smallskip
Let ${\bf R}_t^n$ be the section $\{t\}\times{\bf R}^d$ of ${\bf R}^n$ (in BNF coordinates). 
We look for $K(t,E):L^2({\bf R}^d)\to L^2({\bf R}_t^n)$ (formally), microlocalized near $\Gamma_+(E)$, 
of the form $K(t,E)v(x;h)=\int\int e^{i(S(t,x,\eta)-y\eta)/h}a(t,x,\eta;E,h)v(y)\,dy\wedge d\eta$, and  
such that 
$$H(hD_t,x,hD_x;h)K(t,E)=0, \ K(0,E)=\id$$ 
Considering realizations in the complex domain adapted to the weight $\widetilde\Phi_\theta$, we compute most easily $K(t,E)$ 
in the ``phase of inflation''. Here, 
solving eikonal and transport equations, we find that the leading term of $S$ and $a$ with respect to BNF is given by
those of the Model Hamiltonian, and $K(t,E)$ is also in BNF. Let
$\chi\in C^\infty({\bf R})$, be equal to 0 near 0, 1 near $[2\pi,\infty[$.
There is a $h$-PDO $B(E)=B^w(x,hD_x;E)$ such that 
$L(t,E)=K(t,E)B(E)$ satisfies as in (2.25)
$$\bigl({i\over h}[H,\chi(t)]L(t,E)v|L(t,E)v\bigr)=\bigl(v|v\bigr)\leqno(3.5)$$
Outside the ``phase of inflation'' the analysis is somewhat simpler, since $H-E$ is already elliptic (3.3). 

We set $K_0(t,E)=K(t,E)$ where $K(t,E)$ is Poisson operator with Cauchy data at $t=0$, and $L_0(t,E)=K_0(t,E)B(E)$; we set similarly 
$L_{2\pi}(t,E)=K_0(t-2\pi,E)B(E)$ with Cauchy data at $t=2\pi$. The monodromy operator (or semi-classical Poincar\'e map) is defined by 
$$M^*(E)=L_{2\pi}^*(E){i\over h}[H,\chi]L_0(\cdot,E)\leqno(3.6)$$ 
as an operator on $L^2({\bf R}^d)$, which is a concrete version of (2.28) and (2.33).
As a function de $\chi$, $M^*(E)$ follows a ``0-1 law'': it is 0 if $\supp \chi\subset]0,2\pi[$, and unitary if 
$\chi$ equals 0 near 0, and 1 near $2\pi$. For the model case one has 
$M^*(E)v(x)=e^{-2i\pi E/h}e^{\pi \mu}v(xe^{2\pi \mu})$ since $\int\chi'(t)\, dt=1$.
Unitarity of $M^*(E)$ may not be clear in (3.6), but follows from uniqueness
of the monodromy operator and Proposition 2.1 (when hypotheses match). 
Moreover $M^*(E)$ 
is in BNF, so that eigenfunctions of $M^*(E)$ are
homogeneous polynomials, which leads to Bohr-Sommerfeld quantization rules (see [Lou], [LouRo1,2], 
and a detailed version [LouRo3] in progress). See also [IfaLouRo] for higher order expansions in the 1-D case.
In fact, one can show that 
$M^*(E)=e^{iR^w(x,hD_x;E,h)/h}$,
where $R$ is $h$-PDO in BNF, self-adjoint for real $E$. This gives another proof for unitarity. 
\bigskip
\noindent {\bf Appendix. A short review on complex scaling}
\medskip
Carrying the arguments of [SjZw] to the framework of resonances, 
the proof of Theorem 1.1 in Sect.2 requires only some ``mild'' deformations outside of a neighborhood of $\gamma_0$. 
Sharper deformations are needed in Sect.4 for Theorem 1.2.

For large $x$, the ``dilated'' operator'' takes the form $H_\theta(y,hD_y;h)=U_\theta^*H(y,hD_y;h)U_\theta$.
Here $\theta\in{\bf C}$ is a small parameter ($\im\theta\geq0$  
for outgoing resonances) that we eventually set to $i\theta$ for simplicity).

We say that $U_\theta$ is an {\it analytic dilation} if this is a linear change of variables of the form
$U_\theta u(x)=e^{n\theta/2}u(\theta x)$, and an {\it analytic distorsion} if the change of variables is non linear, but in both cases
it is useful to consider
the scalar product on $L^2({\bf R}^n)$ as a duality product between 
$L^2(\Gamma_\theta)$ and $L^2(\Gamma_{\overline \theta})$ by means of the formula
$$\langle u,v\rangle_\theta=\int_{e^\theta{\bf R}^n}u(y)\overline{v(\overline y)}\,dy\leqno(A.1)$$
For small $\theta\in{\bf C}$, $\Gamma_\theta=e^\theta{\bf R}^n$ is a totally real manifold, whose cotangent space $T^*\Gamma_\theta$,
is a IR-manifold (Lagrangian for $\im d\eta\wedge\,dy$, symplectic for $\re d\eta\wedge\,dy$. 

It makes no difficulty to extend the notion of ``unitary operator'' of ``self-adjoint'' operators in that sense:
for instance if $U_\theta$, for real $\theta$, is unitary on $L^2({\bf R}^n)$, its adjoint for this duality is the analytic extension 
(with respect to small $\theta\in{\bf C}$) of $U_\theta^{-1}$,
and $Q_\theta$ is ``self-adjoint'' means $Q_\theta$
is the analytic continuation of the self-adjoint operator $Q_\theta$ for real $\theta$. 
 
Near $\gamma_0$, $H_\theta$ is defined 
through microlocally weighted $L^2$ (or Sobolev) spaces. 
The microlocal weights $G(y,\eta)$ are chosen among {\it escape functions}, i.e. a smooth functions which is increasing along the flow of $X_H$,
and strictly increasing away from the trapped set; they do not depend, locally, on the energy parameter.
A general result [GeSj] states that there always exists such a function. 
\medskip
\noindent {\it Examples}: (1) Let $H(y,\eta)=\eta^2$, then for any $E>0$, ${\cal K}(E)=\emptyset$,  
and $G(y,\eta)=y\eta$ is an escape function since $X_HG\geq E$
when $|\eta^2-E|\leq E/2$. (2) Let $H(y,\eta)=\eta^2+V(y)$, where $V$ satisfies the virial condition outside a compact set, i.e. 
$2V+y\cdot\nabla V(y)\leq -\delta$ when $y\notin k$. Then $G(y,\eta)$ satisfies $X_HG\geq 2E-2\delta$ when $2|\eta^2+V-E|\leq\delta$.
Modifying it suitably for $y$ close to $k$, so that it vanishes on $k$, we get an escape function outside ${\cal K}(E)=\{(y,\eta): y\in k,
\eta^2+V(y)=E, E>\delta\}$. This is the case (and a paradigm of our situation when restricting to the center manifold)
for $H(y,\eta)=\eta^2-y^2$ where ${\cal K}=\{(0,0)\}$ and $G(y,\eta)=y\eta$. 
\smallskip
In the deformation procedure, escape functions $G(y,\eta)$ have to be modified 
outside a compact set. Namely, for fixed $\lambda>0$, let $G(y,\eta;h)=\lambda h\log(1/h)G_0(y,\eta)$ 
Weighted deformation $h$-PDO $Q(y,\eta;h)$ consists in conjugating
$$Q_G(y,hD_y;h)=e^{-G(y,hD_y;h)/h}Q(y,hD_y;h)e^{G(y,hD_y;h)/h}$$
Due to the mild factor $h\log(1/h)$, $\{Q_G(y,hD_y;h): Q\in S^0(m)\}$ is a ``good'' class of $h$-PDO, bounded on $L^2({\bf R}^n)$.
See [NoSjZw], [NoZw] for details. 

Alternatively (or mixing both techniques) 
complex scaling can be formulated within the theory of $h$-PDO's in the complex domain, where the usual phase space is 
replaced by a IR manifold $\Lambda_\Phi$, and $H(y,hD_y;h)$ is mapped
through a FBI transform to an operator acting on semi-classical distributions microlocalized on $\Lambda_\Phi$.
see [HeSj], [Ma], [Ro].

In Sect.3, we take advantage of BNF to construct escape functions from $G(x,\xi)=x\xi$ in the directions transverse to $\gamma_0$ 
\bigskip

\noindent {\bf References}:

\noindent [A] Marie-Claude Arnaud. On the type of certain periodic orbits minimizing the Lagrangian action. Nonlinearity 11, p.143-150, 1998.

\noindent [BLaz] V.M.Babich, V.Lazutkin. Eigenfunctions concentrated near a closed geodesic. Topics in Math. Phys., Vol.2, M.Birman, ed. 
Consultants' Bureau, New York, p.9-18, 1968. 

\noindent [Bog] E.B.Bogomolny. Semi-classical quantization of multi-dimensional systems. IPNO/TH 91-17, 1991.

\noindent [BrCoDu] P.Briet, J.M.Combes, P.Duclos. On the location of resonances for Schr\"odinger operators II. Comm. Part. Diff. Eq. 12,
p.201-222, 1987.

\noindent [Br] A.D.Bryuno. Normalization of a Hamiltonian system near an invariant cycle or torus. Russian Math. Surveys 44:2, p.53-89, 1991.

\noindent [Chr] H.Christianson. Quantum monodromy and nonconcentration near a closed semi-hyper- bolic orbit. 
Trans. Amer. Math. Soc. 363, No.7, p.3373–3438, 2011. 

\noindent [CdV] Y.Colin de Verdi\`ere. M\'ethodes semi-classiques et th\'eorie spectrale. 
https://www-fourier.ujf-grenoble.fr/~ycolver/ All-Articles/93b.pdf

\noindent [Co] J.M.Combes. Spectral deformations techniques and applications to $N$-body Schr\"odinger operators. Proc. Int. 
Congress of Math. Vancouver, p.369-376, 1974. 

\noindent [FauRoySj] F.Faure, N.Roy, J.Sj\"ostrand. Semi-classical approach for Anosov diffeomorphisms and Ruelle resonances, 2008.

\noindent [GeSi] C.G\'erard, I.M. Sigal. Space-time picture of semi-classical resonances. Comm. Math. Phys. 145, p.281, 1992.
 
\noindent [G\'eSj] C.G\'erard, J.Sj\"ostrand. {\bf 1}. Semiclassical resonances generated by a closed trajectory of hyperbolic type.
Comm. Math. Phys. 108, p.391-421, 1987. {\bf 2}. R\'esonances en limite semiclassique et exposants de Lyapunov.
Comm. Math. Phys. 116, p.193-213, 1988.

\noindent [GuPa] V.Guillemin, T.Paul. Some remarks about semiclassical trace invariants and quantum normal forms. 
Comm. Math. Phys. 294 No. 1, p.1–19, 2010.

\noindent [HeMa] B.Helffer, A.Martinez. Comparaison entre diverses notions de r\'esonances. Helvetica Phys. Acta, 60, p.992-1008, 1987.

\noindent [HeSj] B.Helffer, J.Sj\"ostrand. R\'esonances en limite semi-classique. M\'emoires S.M.F. 114(3), 1986.

\noindent [IfaLouRo] A.Ifa, H.Louati, M.Rouleux. Bohr-Sommerfeld quantization rules revisited: the method of positive commutators,
J. Math. Sci. Univ. Tokyo 25, p.91-137, 2018.

\noindent [Ik] M.Ikawa. Singular perturbation of symbolic flows and poles of the zeta functions. Osaka J.Mats. 27, p.281-300, 1990.

\noindent [Iv] V.Ivrii. Microlocal Analysis and Precise Spectral Asymptotics. Springer-Verlag, Berlin, 1998.

\noindent [KaHa] A.Katok, B.Hasselblatt. Introduction to the modern theorey of dynamical systems. Cambridge Univ. Press, 1999. 

\noindent [Kl] W.Klingenberg. Lectures on closed geodesics. Lect. Notes in Math. 230, Springer.

\noindent [Lou] Hanen Louati. ``R\`egles de quantification semi-classiques pour une orbite p\'eriodique de type hyperbolique''.
Th\`ese, Universit\'es de Toulon et Tunis El-Manar, 2017. 

\noindent [LouRo] H.Louati, M.Rouleux. {\bf 1}. Semi-classical resonances associated with a periodic orbit. 
Math. Notes, Vol. 100, No.5, p.724-730, 2016. {\bf 2}. Semi-classical quantization rules for a periodic orbit of hyperbolic type. 
Proceedings ``Days of Diffraction 2016'', Saint-Petersburg, p.112-117, IEEE.
{\bf 3}. Quantum monodromy and semi-classical quantization rules, {\it in preparation}.

\noindent [Ma] A.Martinez. An introduction to Semiclassical and Microlocal Analysis, Springer, 2001.

\noindent [NoSjZw] S.Nonnenmacher, M.Zworski. Quantum decay rates in chaotic scattering. Acta. Math. 203, p.149-233, 2009.

\noindent [NoSjZw] S.Nonnenmacher, J.Sj\"ostrand, M.Zworski. From Open Quantum Systems to Open Quantum maps.
Comm. Math. Phys. 304, p.1-48, 2011

\noindent [ReSi] M.Reed, B.Simon. Methods of Modern Mathematical Physics IV, Academic Press, 1975.

\noindent [Ro] M.Rouleux. Absence of resonances for semi-classical Schr\"odinger operators of Gevrey type. 
 Hokkaido Math. J., Vol.30 p.475-517, 2001.

\noindent [Sj] J.Sj\"ostrand. {\bf 1}. Singularit\'es analytiques microlocales. Ast\'erisque No.95, 1982. 
{\bf 2}. Resonances associated to a closed hyperbolic trajectory in dimension 2. Asympt. Analysis 36, p.93-113, 2003.
{\bf 3}. Geometric bounds on the density os resonances for semiclassical problems. Duke Math. J. 60, p.1-57, 1990.

\noindent [SjZw] J. Sj\"ostrand and M. Zworski. Quantum monodromy and semi-classical trace formulae. J. Math. Pure Appl.
81(2002), 1-33.

\noindent [Tip] A.Tip. Atoms in circularly polarized fields: 
the dilation-analytic approach. J.Phys. A Math. Gen. Phys. 16, p.3237-3259 (1983)

\noindent [Va] B.Vainberg. On exterior elliptic problems. {\bf I} Mat. Sb. 92(134), 1973, {\bf II} Math. USSR Sb. 21, 1973. 
 
\noindent [Vo] A.Voros. {\bf 1}. Unstable periodic orbits and semiclassical quantization. J.Phys. A(21), p.685-692, 1988.
{\bf 2}. R\'esurgence quantique. Annales Institut Fourier, 43:1509–1534, 1993.
{\bf 3}. Aspects of semiclassical theory in the presence of classical chaos. 
Prog. Theor. Phys. Suppl. No 116, P.17-44, 1994.

\noindent [Zel] S.Zelditch. {\bf 1}. Wave trace invariants at elliptic closed geodesics. GAFA, 7:145–213, 1997.
{\bf 2}.  Wave invariants for non-degenerate closed geodesics. GAFA, 8:179–207, 1998.

\bye